\documentclass[runningheads]{ASMDA}

\usepackage{amsmath,amssymb}
\usepackage{graphicx} 
\def\rmd{\mathrm{d}}

\begin{document}
%
\title*{A strong optimality result for anisotropic self--similar textures}
%
\toctitle{A strong optimality result for anisotropic self--similar textures}
%
\titlerunning{A strong optimality result for anisotropic self--similar textures}
%
\author{
M. Clausel\inst{1}
\and
 B. Vedel\inst{2}
}
%
\index{Clausel,M.}
\index{Vedel,B.}

%
\authorrunning{Clausel et al.}
%
\institute{
Laboratoire Jean Kuntzmann, Universit\'e de Grenoble, CNRS, F38041 Grenoble Cedex 9, France (E-mail: {\tt marianne.clausel@imag.fr})
\and
  Laboratoire de Mathematiques et Applications des Math\'ematiques, Universite de Bretagne Sud, Universit\'e Europ\'eene de Bretagne
Centre Yves Coppens, Bat. B, 1er et., Campus de Tohannic BP 573,
56017 Vannes, France. (E-mail: {\tt vedel@univ-ubs.fr})
  }

\maketitle

\begin{abstract}
In~\cite{RCVJA11,RCVJA13}, we proposed a method to characterize jointly self-similarity and anisotropy properties of a large class of self--similar Gaussian random fields. We provide here a mathematical analysis of our approach,  proving that the sharpest way of measuring smoothness is related to these
anisotropies and thus to the geometry of these fields.
\end{abstract}

\section{Introduction and motivations}
In numerous modern applications (geography \cite{f97}, biomedical imagery (see~\cite{B2010} for example), geophysics \cite{SL87}, art investigation \cite{abry:2012}, \ldots), the data available for analysis consist of images of homogeneous textures, that need to be characterized. For  such images, a key issue consists first  in describing, within a suitable framework, the anisotropy of the texture, and  then in defining regularity anisotropy parameters that can actually and efficiently be measured via numerical procedures and further involved into e.g., classification schemes.

Furthermore, in many cases the analyzed textures  display at the same time self--similarity and anisotropy properties. This is for example such the case in medical imaging (osteoporosis, muscular tissues, mammographies,...), cf. e.g.~\cite{BE03,BMS09}, hydrology~\cite{Pon06}, fracture surfaces analysis~\cite{DH99},\ldots.

In~\cite{RCVJA11,RCVJA13}, we proposed a method for studying jointly selfsimilarity and anisotropy in images by focusing on a specific classic class of Gaussian anisotropic selfsimilar processes.
We consider $\ell^p$ norms of hyperbolic wavelets coefficients which permits the use of different dilation factors along the horizontal and vertical axis. We prove empirically that these $\ell^p$ norms are maximal for a specific ratio of the horizontal and vertical axis, directly related to the anisotropy of the model.

In~\cite{ACJRV13}, we proved that there is a close relationship between $\ell^p$ norms of the hyperbolic wavelet coefficients of a function and its norm in a convenient class of anisotropic functional spaces : the anisotropic Besov spaces. The consistence properties of the estimators introduced in~\cite{RCVJA11,RCVJA13} can then be reformulated using sample paths properties of the model in anisotropic Besov spaces. The aim of this paper is to provide mathematical foundations of the empirical results of~\cite{RCVJA11,RCVJA13}. Here, we prove what we call a {\it strong  optimality result}, namely that the critical exponent of the model in anisotropic Besov spaces is maximal when the parameters of the analyzing space fit these of the analyzed textures, which is the exact transcription into a mathematical way of the empirical results of~\cite{RCVJA11,RCVJA13}.

The paper is organized as follows. In Section~\ref{s:model}, we first present the studied self--similar anisotropic model. Thereafter in Section~\ref{s:besov}, we recall some basic facts about anisotropic Besov spaces. In Section~\ref{s:main} we then state our main result. The proofs are postponed in Section~\ref{s:proof}.

In what follows, we denote $\mathcal{E}^+$ the collection of $2\times 2$ matrices whose eigenvalues have positive real parts.
\section{Presentation of the studied model}\label{s:model}
In~\cite{RCVJA11,RCVJA13}, we choosed to investigate the properties of a large class of anisotropic Gaussian self-similar fields, introduced in~\cite{BE03,BMS09}, referred as \emph{Operator Scaling Gaussian Random Field}, in short OSGRF.

For any matrix $E_0$ belonging in $\mathcal{E}^+$ such that $\mathrm{Tr}(E_0)=2$ and any $H_0\in (0,\min_{\lambda\in Sp(E_0)}\mathrm{Re}(\lambda))$, this class can be defined using the following harmonizable representation:
\begin{equation}\label{d:model}
X_{\rho,E_0,H_0}(\underline{x}) =  \int_{\mathbb{R}^2} (\mathrm{e}^{\mathrm{i}\langle \underline{x}, \, \underline{\xi} \rangle} -1) \rho(\underline{\xi})^{-(H_0 +1)} d\widehat{W}(\underline{\xi})\;,
\end{equation}
where $\underline{x}=(x_1,x_2)$,  $\underline{\xi}=(\xi_1,\xi_2)$ and $d\widehat{W}(\underline{\xi})$ stands for a 2D  {Wiener measure}. Here $\rho$ is an $E_0$--homogeneous continuous positive function,that is satisfies the following homogeneity relationship
$\rho(a^{E_0}\xi)=af(\xi)$ on $\mathbb{R}^2$. We also assume that $\int (1\wedge |\xi|^2)\rho(\xi)^{-2(H_0 +1)} d\xi<+\infty$ which ensures the existence of the Gaussian field $X$ (see~\cite{BMS09}).

With this construction, the OSGRF $X_{\rho,E_0,H_0}$ has stationary increments. Further it satisfies an anisotropic version of the scaling property with parameter $H_0$ (where $\overset{\mathcal{L}}{=}  $ denotes equality for all finite dimensional distributions):
\begin{equation}\label{e:selfsimilar}
\forall a>0,\,\{ X_{\rho}(a^{E_0}  \underline{x} )\}\overset{\mathcal{L}}{=} \{ a^{H_0} X_{\rho}( \underline{x}) \}\;.
\end{equation}
(as usual $a^{E_0}=\exp(E_0\log(a))=\sum_k (\log a)^k E_0^k/k!$).

Constructions of $E_0$--homogeneous continuous positive function have been proposed in~\cite{BMS09} via an integral formula (Theorem~2.11). An alternative construction, more fitted for numerical simulations, can be found in~\cite{CV11}.

Let us now give a concrete example of OSRGF. Consider the case where $E_0=\begin{pmatrix}\alpha_0&0\\0&2-\alpha_0\end{pmatrix}$ with $\alpha_0\in (0,2)$ and set
\[
\rho(\xi_1,\xi_2)=|\xi_1|^{1/\alpha_0}+|\xi_2|^{1/(2-\alpha_0)}\;.
\]
The function $\rho$ is obviously an $E_0$ homogeneous positive continuous fonction. For any $H_0\in \min(\alpha_0,2-\alpha_0)$, the associated Gaussian field will be
\[
X_{\rho,E_0,H_0}(x_1,x_2)=\int_{\mathbb{R}^2}\frac{\mathrm{e}^{\mathrm{i}\langle \underline{x}, \, \underline{\xi} \rangle} -1}{(|\xi_1|^{1/\alpha_0}+|\xi_2|^{1/(2-\alpha_0)})^{H_0+1}}d\widehat{W}(\underline{\xi})\;,
\]
The scaling property satisfied by this field is then 
\[
\forall a>0,\,\{ X_{\rho,E_0,H_0}(a^{\alpha_0}x_1,a^{2-\alpha_0}x_2)\}\overset{\mathcal{L}}{=} \{ a^{H_0} X_{\rho,E_0,H_0}(x_1,x_2) \}\;.
\]
\section{Anisotropic concepts of smoothness}\label{s:besov}
Our main goal here is to study the sample paths properties of this class of
Gaussian fields in suitable anisotropic functional spaces. This
approach is quite natural (see~\cite{Kam96}) since the studied model is anisotropic. To this end, suitable concepts of anisotropic smoothness are needed. The aim of this section is to give some background about the appropriate anisotropic functional spaces : Anisotropic
Besov spaces. These spaces generalize classical (isotropic) Besov
spaces and have been studied in parallel with them.
Let $D$ a diagonalizable matrix of $\mathcal{E}^+$ with eigenvalues $\lambda_1,\lambda_2$ and associated eigenvalues $(e_1,e_2)$. The anisotropic Besov spaces with anisotropy $D$ can be defined as follows (see Theorem~5.8 of \cite{Trie06})~:
\begin{definition}\label{d:besov}
Let $(p,q)\in [1,+\infty]^2$, $s>0$ and $\beta\in \mathbb{R}$. For any $f\in\L^p(\mathbb{R}^2)$ set
\[
\|f\|_{\dot{B}^{s}_{p,q,|\log|^\beta}(\mathbb{R}^{2},D)}=\sum_{\ell=1}^2\left(\int_{0}^1 \|(\Delta^{M_\ell}_{t e_\ell} f)(x)\|_{L^p}^q
t^{-sq/\alpha_{\ell}-1}|\log(t)|^{-\beta q/\alpha_{\ell}}\rmd t\right)^{1/q}\;.
\]
By definition
\[
B^{s}_{p,q,|\log|^\beta}(\mathbb{R}^{2},D)=\{f\in\L^p(\mathbb{R}^2),
\|f\|_{\dot{B}^{s}_{p,q,|\log|^\beta}(\mathbb{R}^{2},E)}<+\infty\}\;.
\]
The matrix $D$ is called the anisotropy of the Besov space $B^{s}_{p,q,|\log|^\beta}(\mathbb{R}^{2},D)$.
\end{definition}
\begin{remark}
It is well--known that there exists a strong relationship between $\ell^p$ norms and classical wavelet coefficients (see~\cite{Trie83,Trie06}). In~\cite{RCVJA11,RCVJA13}, our estimators are based on  hyperbolic wavelet analysis. In~\cite{ACJRV13}, we proved that $\ell^p$ norms of hyperbolic wavelet coefficients of a functions are related to its norms in anisotropic Besov spaces. It is the reason why the natural mathematical framework to relate anisotropy and self--similar properties of the model to its sample paths properties is this of anisotropic Besov spaces.
\end{remark}
\begin{remark}\label{rem:comp} Let $D$ a diagonalizable matrix of $\mathcal{E}^+$. For any
$a>0$, $\lambda D$ is also a diagonalizable matrix of $\mathcal{E}^+$ with eigenvalues $a\lambda_1, a\lambda_2$ and same eigenvectors as $D$. Hence for any $s>0$,
$B^{a s}_{p,q}(\mathbb{R}^2,\lambda D)=B^{s}_{p,q}(\mathbb{R}^2, D)$.
\end{remark}
Using Remark~\ref{rem:comp}, we deduce that without loss of generality, we can assume in the sequel that
$\mathrm{Tr}(D)=2$. We then define
\[
\mathcal{E}^+_2=\{E\in\mathcal{E}^+,\,\mathrm{Tr}(E)=2\}\;.
\]
As it is the case for isotropic spaces, anisotropic H\"older
spaces $\mathcal{C}^s(\mathbb{R}^2, E)$ can be defined as particular anisotropic Besov
spaces.
\begin{definition}
Let $s>0$, $\beta\in\mathbb{R}$. The anisotropic
H\"older spaces ${\mathcal{C}}^s_{|\log|^\beta}(\mathbb{R}^2, D)$ are defined by
\[
\mathcal{C}^{s}_{|\log|^{\beta}}(\mathbb{R}^{2},E)=B^{s}_{\infty,\infty,|\log|^{\beta}}(\mathbb{R}^{2},E)\;.
\]
\end{definition}
Hence, a bounded function $f$
belongs to $\mathcal{C}^{s}_{|\log|^{\beta}}(\mathbb{R}^{2},E)$ if and
only if for any $r\in (0,1)$, $\Theta\in S_0^{E}(|\cdot|_E)$ and
$x\in \mathbb{R}^d$, one has
\[
|f(x+r^E \Theta)-f(x)|\leq C_0 r^s|\log(r)|^\beta\;,
\]
for some $C_0>0$, that is if and only if its
restriction $f_{\Theta}$ along any parametric curve of the form
\[
r>0\mapsto r^E \Theta\;,
\]
with $\Theta\in S_0^{E}(|\cdot|_E)$ is in the usual H\"{o}lder
space $\mathcal{C}^{s}_{\vert \log \vert^\beta}(\mathbb{R})$ and
$\|f_{\Theta}\|_{\mathcal{C}^{s}_{\vert \log \vert^\beta}(\mathbb{R})}$
does not depend on $\Theta$. Roughly speaking, the anisotropic
``directional'' regularity in any anisotropic ``direction'' has to
be larger than $s$. In other words, we replace straight lines of
isotropic setting by curves with parametric equation $r>0\mapsto r^E \Theta$ adapted to anisotropic setting.

To state our optimality results we need a local version of
anisotropic Besov spaces~:
\begin{definition}
Let $D \in \mathcal{E}^{+}$ be a fixed diagonalizable anisotropy, $1\leq p,q\leq \infty$, $\beta\in\mathbb{R}$, $0<s<\infty$ and $f\in L^{p}_{loc}(\mathbb{R}^{2})$.

The function $f$ belongs to
$B^{\alpha}_{p,q,|\log|^\beta,loc}(\mathbb{R}^{2},E)$ if for any $\varphi\in
\mathcal{D}(\mathbb{R}^{d})$, the
function $\varphi f$ belongs to $B^{\alpha}_{p,q,|\log|^\beta}(\mathbb{R}^{2},D)$.

The anisotropic local critical exponent in anisotropic Besov
spaces $B^{s}_{p,q}(\mathbb{R}^2,E)$ of $f\in
L^p_{loc}(\mathbb{R}^{2})$ is then defined by
\[
\alpha_{f,loc}(E,p,q)=\sup\{s,\,f\in
B^{s}_{p,q,loc}(\mathbb{R}^2,E)\}\;.
\]
\end{definition}
\section{Statement of our main result}\label{s:main}
In what follows, we are given $E_0\in\mathcal{E}^+_2$ and $\rho_{E_0}$
an $E_0^{\!\!\!\!\!\!\! t}\;$--homogeneous continuous positive function, $H_0\in (0,\min_{\lambda\in Sp(E_0)}(\mathrm{Re}(\lambda))$. 

Our results will be based on a comparison
between the topology related to $\rho_{E_0}$
involved in the construction of the Gaussian field
$\{X_{\rho_{E_0},E_0,H_0}(x)\}_{x\in\mathbb{R}^d}$ defined by equation~(\ref{d:model}) and this of the
analyzing spaces $B^s_{p,q}(\mathbb{R}^2,D)$. To be able to
compare these two topologies, we also assume that
$D\in\mathcal{E}^+_2$.

We characterize in some sense an anisotropy $E_0$ and an Hurst index of the field $\{X_{\rho_{E_0},H_0}(x)\}_{x\in\mathbb{R}^d}$, which is interesting when analyzing anisotropic self--similar textures. Combining the results of this paper and these of~\cite{ACJRV13}, our approach can thus be turned in an effective algorithm for the estimation of the anisotropy of self similar textures (see~\cite{RCVJA11,RCVJA13}). In the case where $q=\infty$, Theorem~1 is the exact mathematical reformulation of equation (9) stated in~\cite{RCVJA13}, on which are based the definition of the estimators of the anisotropy and smoothness of the model. Our main result in then both a mathematical justification and an extension of the results empirically proved in~\cite{RCVJA11,RCVJA13}~:
\begin{theorem}\label{ThOptim} Let $(p,q)\in [1,+\infty]^2$. Then almost surely
\[
\alpha_{X_{\rho_{E_0},E_0,H_0},loc}(E_0,p,q)=\sup\{\alpha_{X_{\rho_{E_0},E_0,H_0},loc}(D,p,q),\,D\mbox{ diagonalizable and }D\in
\mathcal{E}^{+}_2\}=H_0\;.
\]
\end{theorem}
In fact, Theorem~\ref{ThOptim} contains two
main results :
\begin{itemize}
\item The critical exponent of the field $\{X_{\rho_{E_0},H_0}(x)\}_{x\in\mathbb{R}^{d}}$ in
anisotropic Besov space $B^{s}_{p,q}(\mathbb{R}^d,D_0)$ where $D_0$ is the real diagonalizable part of $E_0$ (see Proposition 4.1 of~\cite{CV13b} for a definition) equals the
associated Hurst index $H_0$. It has already been proved in \cite{CV13b} in a general $d$--dimensional context.
\item The diagonalizable real part $D_0$ of any anisotropy $E_0$ of the field $\{X_{\rho_{E_0},H_0}(x)\}_{x\in\mathbb{R}^{d}}$ maximizes
this critical exponent among all possible analysis
matrices. In fact, the ``best way'' of measuring smoothness of the
field $\{X_{\rho_{E_0},H_0}(x)\}_{x\in\mathbb{R}^{d}}$ is to measure smoothness along the ``anisotropic
directions'' $r>0\mapsto r^{D_0}\Theta$, related to the genuine geometry of the field.
\end{itemize}
\section{Proof of Theorem~\ref{ThOptim}}\label{s:proof}
Proposition~5.2 of \cite{CV13b} directly implies that a.s.
\[
\alpha_{X_{\rho_{E_0},H_0},loc}(E_0,p,q)=H_0\;.
\]
We now prove that if $D$ is a diagonalizable matrix of $\mathcal{E}^+_2$ then a.s.
\[
\alpha_{X_{\rho_{E_0},H_0},loc}(D,p,q)\leq H_0\;.
\]
Denote $D_0$ the real diagonalizable part of $E_0$ (see Proposition 4.1 of~\cite{CV13b}). Observe now that, by
Theorem 1.1 of \cite{CV13a} and Lemma~5.1 of \cite{CV13b}, one has
a.s.
\[
\alpha_{X_{\rho_{E_0},H_0}}(D,p,q)=\alpha_{X_{\rho_{D_0},H_0}}(D,p,q)\;,
\]
where $\rho_{E_0}$, $\rho_{D_0}$ are respectively two
$E_0$-- and $D_0$--homogeneous continuous positive functions
and $\{X_{\rho_{E_0},H_0}(x)\}_{x\in\mathbb{R}^{d}}$, $\{X_{\rho_{D_0},H_0}(x)\}_{x\in\mathbb{R}^{d}}$ the associated
Gaussian fields by~(\ref{d:model}). One can then investigate the
sample paths properties of
$\{X_{\rho_{D_0},H_0}(x)\}_{x\in\mathbb{R}^{d}}$ instead of this
of $\{X_{\rho_{E_0},H_0}(x)\}_{x\in\mathbb{R}^{d}}$.

Hence from now, we then assume that $E_0$ equals its
diagonalizable real part, namely that $E_0=D_0$. Even if $D$ and $D_0$ are both diagonalizable since these matrices are not commuting we cannot assume
that $D_0$ and $D$ are both diagonal.

Let $\varphi\in \mathcal{D}(\mathbb{R}^{d})$. As above, one may
assume $\mathrm{supp}(\varphi)\subset K=B_{D_0}(0,1)$. Denote
$\lambda_1^0\leq \lambda_2^0$ (resp $\lambda_1\leq \lambda_2$) the
two eigenvalues of $D_0$ (resp. $D$) and $e_1^0,e_2^0$ (resp
$e_1,e_2$) some associated eigenvectors. We exclude the two cases $\lambda_1^0=\lambda_2^0$,
$\lambda_1=\lambda_2$ corresponding to the cases $D_0=Id$, $D=Id$ which can
be deduced from Theorem 4.1 of \cite{CV13b}. Hence the inequalities between the eigenvalues are strict.

Let us fix $p\in
[1,+\infty]$. By assumption $[H_0/\lambda_1^0]=[H_0/\lambda_2^0]=1$. The proof of
Proposition~5.5 of \cite{CV13b} provides us with the following equalities
\begin{equation}\label{e:h0}
\begin{array}{lll}
H_0&=&\sup\{\beta, \int_{0}^1\|(\varphi X)(x+t
e_2^0)-(\varphi X)(x)\|_{L^p}^p t^{-\beta p/\lambda_1^0-1}\mathrm{d}t<+\infty\}\\
&=&\sup\{\beta, \int_{0}^1\|(\varphi X)(x+t
e_1^0)-(\varphi X)(x)\|_{L^p}^p t^{-\beta p/\lambda_2^0-1}\mathrm{d}t<+\infty\}\;.
\end{array}
\end{equation}

Define now the two followings indices
\[
\alpha_1(p)=\sup\{\beta, \int_{0}^1\|(\varphi X)(x+t e_1)-(\varphi
X)(x)\|_{L^p}^p t^{-\beta p/\lambda_1-1}\mathrm{d}t<+\infty\}\;,
\]
\[
\alpha_2(p)=\sup\{\beta, \int_{0}^1\|(\varphi X)(x+t e_2)-(\varphi
X)(x)\|_{L^p}^p t^{-\beta p/\lambda_2-1}\mathrm{d}t<+\infty\}\;.
\]

Thereafter set
$\alpha(p)=\min(\alpha_1(p),\alpha_2(p))$. Assume that for some
$p\in [1,+\infty]$ one has a.s. $\alpha(p)>H_0$, that is a.s $\varphi X\in
B^{\alpha(p)}_{p,p}(\mathbb{R}^{2},D)$ with $\alpha(p)>H_0$.

We first need the following lemma~:
\begin{lemma}\label{lem:sep}
Only five cases are possible~:
\begin{enumerate}
\item
\[
\frac{H_0}{\lambda_1^0}=\max\left(\frac{H_0}{\lambda_1^0},\frac{H_0}{\lambda_2^0}\right)\leq \min\left(\frac{\alpha_1(p)}{\lambda_1},\frac{\alpha_2(p)}{\lambda_2}\right)\;.
\]
\item
\[
\frac{H_0}{\lambda_1^0}=\max\left(\frac{H_0}{\lambda_1^0},\frac{H_0}{\lambda_2^0}\right)>
\min\left(\frac{\alpha_1(p)}{\lambda_1},\frac{\alpha_2(p)}{\lambda_2}\right)>\frac{H_0}{\lambda_2^0}=\min\left(\frac{H_0}{\lambda_1^0},\frac{H_0}{\lambda_2^0}\right) \;.
\]
\item
\[
\frac{H_0}{\lambda_1^0}=\max\left(\frac{H_0}{\lambda_1^0},\frac{H_0}{\lambda_2^0}\right)>
\frac{H_0}{\lambda_2^0}=\min\left(\frac{H_0}{\lambda_1^0},\frac{H_0}{\lambda_2^0}\right) \min\left(\frac{\alpha_1(p)}{\lambda_1},\frac{\alpha_2(p)}{\lambda_2}\right)\;.
\]
\item
\[
\frac{H_0}{\lambda_1^0}=\max\left(\frac{H_0}{\lambda_1^0},\frac{H_0}{\lambda_2^0}\right)>\max\left(\frac{\alpha_1(p)}{\lambda_1},\frac{\alpha_2(p)}{\lambda_2}\right)>
\min\left(\frac{\alpha_1(p)}{\lambda_1},\frac{\alpha_2(p)}{\lambda_2}\right)=\frac{H_0}{\lambda_2^0}=\min\left(\frac{H_0}{\lambda_1^0},\frac{H_0}{\lambda_2^0}\right) \;.
\]
\item
\[
\max\left(\frac{\alpha_1(p)}{\lambda_1},\frac{\alpha_2(p)}{\lambda_2}\right)>\frac{H_0}{\lambda_1^0}=\max\left(\frac{H_0}{\lambda_1^0},\frac{H_0}{\lambda_2^0}\right)>
\min\left(\frac{\alpha_1(p)}{\lambda_1},\frac{\alpha_2(p)}{\lambda_2}\right)=\frac{H_0}{\lambda_2^0}=\min\left(\frac{H_0}{\lambda_1^0},\frac{H_0}{\lambda_2^0}\right) \;.
\]
\end{enumerate}
\end{lemma}
{\bf Proof of Lemma~\ref{lem:sep}.} The only point to prove is that the case
\[
\max(\frac{\alpha_1(p)}{\lambda_1},\frac{\alpha_2(p)}{\lambda_2})=\frac{H_0}{\lambda_1^0}=\max(\frac{H_0}{\lambda_1^0},\frac{H_0}{\lambda_2^0})>
\min(\frac{\alpha_1(p)}{\lambda_1},\frac{\alpha_2(p)}{\lambda_2})=\frac{H_0}{\lambda_2^0}=\min(\frac{H_0}{\lambda_1^0},\frac{H_0}{\lambda_2^0}) \;.
\]
is impossible. Suppose that this relationship holds. Then one has
\[
\left(\frac{\alpha_1(p)}{\lambda_1}=\frac{H_0}{\lambda_1^0}\mbox{  and  }\frac{\alpha_2(p)}{\lambda_2}=\frac{H_0}{\lambda_2^0}\right)\mbox{ or }\left(\frac{\alpha_2(p)}{\lambda_2}=\frac{H_0}{\lambda_1^0}\mbox{  and  }\frac{\alpha_1(p)}{\lambda_1}=\frac{H_0}{\lambda_2^0}\right)\;.
\]
Since $\min(\alpha_1(p),\alpha_2(p)>H_0$, it implies that
\[
\left(\frac{H_0}{\lambda_1}<\frac{H_0}{\lambda_1^0}\mbox{  and  }\frac{H_0}{\lambda_2}<\frac{H_0}{\lambda_2^0}\right)\mbox{ or }\left(\frac{H_0}{\lambda_2}<\frac{H_0}{\lambda_1^0}\mbox{  and  }\frac{H_0}{\lambda_1}<\frac{H_0}{\lambda_2^0}\;,
\right)
\]
that is
\[
\left(\lambda_1>\lambda_1^0\mbox{  and  }\lambda_2>\lambda_2^0\right)\mbox{ or }\left(\lambda_2>\lambda_1^0\mbox{  and  }\lambda_1>\lambda_2^0\right)\;.
\]
Since $\lambda_1^0+\lambda_2^0=\lambda_1+\lambda_2=2$, it yields to a contradiction.

We now deal successively with the five cases of Lemma~\ref{lem:sep}. Let us first assume that point~(1) of Lemma~\ref{lem:sep} holds.

Let $s\in (H_0/\lambda_1^0,\min(\alpha_1(p)/\lambda_1,\alpha_2(p)/\lambda_2,1))$. Then, by definition of $s$ one has
\[
\int_{0}^1\|(\varphi X)(x+t e_1)-(\varphi
X)(x)\|_{L^p}^p t^{-s p-1}\mathrm{d}t<+\infty\;,
\]
and
\[
\int_{0}^1\|(\varphi X)(x+t e_2)-(\varphi
X)(x)\|_{L^p}^p t^{-s p-1}\mathrm{d}t<+\infty\;.
\]
The, using the finite difference characterization of classical Besov spaces, one deduces that $\varphi X\in B^s_{p,p}(\mathbb{R}^2)$. Since, by Proposition~5.6 of~\cite{CV13b}, $B^s_{p,p}(\mathbb{R}^2)\hookrightarrow B^{s\lambda_1^0}_{p,p}(\mathbb{R}^2,D_0)$ and $s>H_0/\lambda_1^0$ it yields to a contradiction. \\

\noindent Now suppose that point~(2) of Lemma~\ref{lem:sep} holds. Assume that $\min(\alpha_1(p)/\lambda_1,\alpha_2(p)/\lambda_2)=\alpha_1(p)/\lambda_1$. The other case will be similar.  
Observe that that if $e_1^0=a e_1$
for some $a\in\mathbb{R}$, one then has
\[
H_0/\lambda_1^0=\sup\{s, \int_{0}^1\|(\varphi X)(x+t e_1)-(\varphi
X)(x)\|_{L^p}^p
t^{-sp-1}\mathrm{d}t<+\infty\}=\alpha_1(p)/\lambda_1\;,
\]
which is impossible since by assumption
$H_0/\lambda_1^0>\alpha_1(p)/\lambda_1$. 

Then the family $(e_1,e_1^0)$ is necessarily a basis of
$\mathbb{R}^2$. Hence 
\begin{equation}\label{e:e02}
e_2^0=a e_1+b e_1^0\;,
\end{equation} 
for some $(a,b)\in
\mathbb{R}^{2}$. By ~(\ref{e:h0}), for any
$\varepsilon>0$ one has a.s.
\begin{equation}\label{e:Xloc1}
\int_{0}^1\|(\varphi X)(x+t e_2^0)-(\varphi X)(x)\|_{L^p}^p
t^{-(H_0 p/\lambda_2^0+\varepsilon)-1}\mathrm{d}t=+\infty\;.
\end{equation}
Further, the triangular inequality and equation~(\ref{e:e02}) implies that
\begin{equation}\begin{array}{lll}\label{e:Xtriang}
&&\int_{0}^1\|(\varphi X)(x+t e_2^0)-(\varphi X)(x)\|_{L^p}^p
t^{-(H_0 p/\lambda_2^0+\varepsilon)-1}\mathrm{d}t\\
&\leq& 2^{p-1}\int_{0}^1\|(\varphi X)(x+t(a e_1+b e_1^0))-(\varphi X)(x+t
b e_1^0)\|_{L^p}^p
t^{-(H_0 p/\lambda_2^0+\varepsilon)-1}\mathrm{d}t\\
&&+2^{p-1}\int_{0}^1\|(\varphi X)(x+t b e_1^0)-(\varphi X)(x)\|_{L^p}^p
t^{-(H_0 p/\lambda_2^0+\varepsilon)-1}\mathrm{d}t\;.
\end{array}
\end{equation}
Set now $y=x+t b e_1^0$ and remark that
\[
\|(\varphi X)(x+t(a e_1+b e_1^0))-(\varphi X)(x+t b
e_1^0)\|_{L^p}^p=\|(\varphi X)(y+t a e_1)-(\varphi X)(y)\|_{L^p}^p
\]
In addition, $H_0/\lambda_2^0<\alpha_1/\lambda_1<H_0/\lambda_0^1$.  Since
\[
H_0/\lambda_0^1=\sup\{s, \int_{0}^1\|(\varphi X)(x+t e_1)-(\varphi
X)(x)\|_{L^p}^p t^{-sp-1}\mathrm{d}t<+\infty\}\;,
\] and by definition of $\alpha_1(p)$, one deduces that for any $\varepsilon>0$
sufficiently small
\begin{equation}\label{e:Xloc2}
\int_{0}^1\|(\varphi X)(y+t a e_1)-(\varphi
X)(y)\|_{L^p}^pt^{-(\frac{H_0 p}{\lambda_2^0}+\varepsilon)-1}\mathrm{d}t\leq
\int_{0}^1\|(\varphi X)(y+t a e_1)-(\varphi
X)(y)\|_{L^p}^pt^{-(\frac{\alpha_1(p)p}{\lambda_1}-\varepsilon)-1}\mathrm{d}t\;,
\end{equation}
and
\begin{equation}\label{e:Xloc3}
\int_{0}^1\|(\varphi X)(x+t b e_1^0)-(\varphi X)(x)\|_{L^p}^p
t^{-(\frac{H_0 p}{\lambda_2^0}+\varepsilon)-1}\mathrm{d}t\leq \int_{0}^1\|(\varphi X)(x+t b e_1^0)-(\varphi X)(x)\|_{L^p}^p
t^{-(\frac{H_0 p}{\lambda_1^0}-\varepsilon)-1}\mathrm{d}t\;.
\end{equation}
are finite. Combining this information with Equations~(\ref{e:Xtriang}), (\ref{e:Xloc1}), implies a contradiction.
The proof of Theorem~\ref{ThOptim} in the three other cases of Lemma~\ref{lem:sep} is exactly similar.


{\footnotesize
\begin{thebibliography}{100}
\bibitem{abry:2012} {\sc P. Abry et al.} (2012). When {V}an {G}ogh meets {M}andelbrot: Multifractal Classification of Painting Textures, {\em Signal Processing} To appear.
\bibitem{ACJRV13} {\sc P. Abry, M. Clausel, S. Jaffard, S.G. Roux and B. Vedel} (2013). Hyperbolic wavelet transform: an efficient tool for multifractal analysis of anisotropic textures. {\em Submitted}.
\bibitem{B2010} {\sc M.Bergounioux and  L. Piffet} (2010). A second-order model for image denoising, {\em Set Valued and Variational Analysis},18 (3--4),277--306.
\bibitem{BE03} {\sc A. Bonami and A. Estrade} (2003). Anisotropic analysis of some Gaussian models. {\em The Journal of Fourier Analysis and Applications}  9, 215-236.
\bibitem{BMS09} {\sc H. Bierm\'e, M.M. Meerschaert and H.P. Scheffler} (2009). Operator Scaling Stable
Random Fields. {\em Stoch. Proc. Appl.}  117(3) 312--332.
\bibitem{CV11} {\sc M. Clausel and B. Vedel} (2011). Explicit constructions of operator scaling self--similar random Gaussian fields
{\em Fractals} 19(1) 101--111.
\bibitem{CV13a} {\sc M. Clausel and B. Vedel} (2013). Sample paths properties of Gaussian fields with equivalent spectral densities.
{\em Submitted}.
\bibitem{CV13b} {\sc M. Clausel and B. Vedel} (2013). An optimality result about sample path properties of Operator Scaling Gaussian Random Fields
{\em Submitted}.
\bibitem{DH99} {\sc S. Davies and P. Hall} (1999). Fractal analysis of surface
roughness by using spatial data (with discussion). {\em J. Roy. Statist. Soc. Ser.}  B 61 3--37.
\bibitem{f97} {\sc P. Frankhauser} (1997). L'approche fractale : un nouvel outil dans l'analyse spatiale des
	agglomerations urbaines. {\em Population}, 4 1005--1040.
\bibitem{Kam96} {\sc A. Kamont} (1996). On the Fractional Anisotropic
Wiener Field. {\em Prob. and Math. Stat.}  16(1) 85--98.
\bibitem{Pon06} {\sc L. Ponson et al.} (2006). Anisotropic self-affine properties of experimental fracture surfaces. {\em Int.Journ. of fracture} 140 27--37.
\bibitem{RCVJA11} {\sc S.G. Roux, M. Clausel, B. Vedel, S. Jaffard and P. Abry} (2011). Transform\'ee en ondelettes hyperboliques pour la caract\'erisation des images autosimilaires anisotropes {\em XXIII ieme Colloque GRETSI, Bordeaux 2011}.
\bibitem{RCVJA13} {\sc S.G. Roux, M. Clausel, B. Vedel, S. Jaffard and P. Abry} (2013). The {H}yperbolic
{W}avelet {T}ransform for self-similar anisotropic texture analysis {\em Submitted}.
\bibitem{SL87} {\sc D. Schertzer and S. Lovejoy} (1987). Physically based rain and cloud modeling by anisotropic, multiplicative turbulent cascades.{\em J. Geophys. Res.}. 92 9693--9714.
\bibitem{Trie83} {\sc H. Triebel} (1983). Theory of Function Spaces. Monographs in Math. 78
{\em Birkh\"{a}user Verlag}.
\bibitem{Trie06} {\sc H. Triebel} (2006). Theory of functions spaces
 III {\em Birkha\"{u}ser Verlag}.
\end{thebibliography}
}
\end{document}